\documentclass[12pt]{article}
\title{ On  $\psi$-Appell polynomials and $Q(\partial_{\psi})$-difference calculus nonhomogeneous equation}
\author{A.K.Kwa\'sniewski \\
\\Higher School of Mathematicand Applied Informatics\\
PL-15-021 Bia{\l}ystok, ul.Kamienna 17, POLAND\\
e-mail: kwandr@uwb.edu.pl}

\usepackage{amsmath,amsthm}

\chardef\bslash=`\\ 
\newtheorem{prop}{Proposition}[section]

\begin{document}
\maketitle
\begin{abstract}
One discovers why Morgan Ward solution \cite{1} of $\psi$- {\em
difference calculus} nonhomogeneous equation
$\Delta_{\psi}f=\varphi$ in the form $$f(x)=\sum_{n \geq
1}\frac{B_{n}}{n_{\psi}!}\varphi^{(n-1)}(x)+\int_{\psi}\varphi(x)
+p(x)$$
 recently proposed by the present author (see-below) - extends here now to
$\psi$- Appell polynomials case - almost {\em automatically}. The
reason for that is just proper framework i.e. that of the
$\psi$-{\em Extended Finite Operator Calculus}\textbf{(EFOC)
}recently being developed and promoted by the present author
\cite{2,3,4,5}. Illustrative specifications to $q$-calculus case
and Fibonomial calculus case \cite{5,6} were already made explicit
in \cite{8} due to the of upside down convenient notation for
objects of \textbf{EFOC }as to be compared with functional
formulation [9].

\end{abstract}
A.M.S Classification numbers: 11B39, 11B65, 05A15

\section{Remark on references usage and the "upside down" notation .}
At first let us make a remark on notation.$\psi$ is a number or
functions`  sequence - sequence of functions of a parameter $q$.
$\psi$ denotes an extension of $\langle\frac{1}{n!}\rangle_{n\geq
0}$ sequence to quite arbitrary one (the so called - "admissible"
[2-5]).The specific choices are for example : Fibonomialy-extended
sequence $\langle\frac{1}{F_n!}\rangle_{n\geq 0}$ ($\langle F_n
\rangle$ - Fibonacci sequence )  or just "the usual"
$\psi$-sequence $\langle\frac{1}{n!}\rangle_{n\geq 0}$ or Gauss
$q$-extended $\langle\frac{1}{n_q!}\rangle_{n\geq 0}$ admissible
sequence of extended umbral operator calculus, where $
n_q=\frac{1-q^n}{1-q}$ and $n_q!=n_q(n-1)_q! , 0_q!=1$.

The simplicity of calculations is being achieved due to writing
objects of these extensions in mnemonic convenient \textbf{upside
down notation} \cite{2} , \cite{5}
\begin{equation}\label{}
\frac {\psi_{(n-1)}}{\psi_n}\equiv n_\psi,
n_\psi!=n_\psi(n-1)_\psi!, n>0 ,   x_{\psi}\equiv \frac
{\psi{(x-1)}}{\psi(x)} ,
\end{equation}
\begin{equation}\label{}
x_{\psi}^{\underline{k}}=x_{\psi}(x-1)_\psi(x-2)_{\psi}...(x-k+1)_{\psi}
\end{equation}
\begin{equation}\label{}
x_{\psi}(x-1)_{\psi}...(x-k+1)_{\psi}=
\frac{\psi(x-1)\psi(x-2)...\psi(x-k)} {\psi(x)
\psi(x-1)...\psi(x-k +1)} .
\end{equation}

If one writes the above in the form $x_{\psi} \equiv \frac
{\psi{(x-1)}}{\psi(x)}\equiv \Phi(x)\equiv\Phi_x\equiv x_{\Phi}$ ,
one sees that the name upside down notation is legitimate.

As for references - the papers of main references are:
\cite{1,2,3}. Consequently we shall then take here notation from
\cite{2,3} and the results from \cite{1} as well as from
\cite{2,3}- for granted. For other references see: \cite{2,3,5}
(Note the access via ArXiv to \cite{3,5}).

$A_{n}$ denotes here $\psi$- {\em Appell-Ward numbers} -
introduced below.
\section{ $Q(\partial_{\psi})$- difference  nonhomogeneous equation}
Let us recall \cite{3,2} the simple fact.
\begin{prop}
$Q(\partial_{\psi})$ is a $\psi$- delta operator iff there exists
invertible $S\in \Sigma _{\psi}$ such that
$Q(\partial_{\psi})=\partial_{\psi}S$.
\end{prop}
Formally: "$S=Q/\partial_{\psi}$" or "$S^{-1}=\partial_{\psi}/Q$".
In the sequel we use this abbreviation $Q(\partial_{\psi})\equiv
Q$.

$\psi$- {\em Appell} or generalized {\em Appell polynomials}
$\left\{A_{n}(x)\right\}_{n \geq 0}$ are defined according to
\begin{equation} \label{psi-A-1}
\partial_{\psi}A_{n}(x)=n_{\psi}A_{n-1}(x)
\end{equation}
and they naturally do satisfy the $\psi$ - Sheffer-Appell identity
\cite{3,2}
\begin{equation}\label{psi-A-2}
A_{n}(x+_{\psi}y)=\sum_{s=0}^{n}\binom{n}{s}
_{\psi}A_{s}(y)x^{n-s}.
\end{equation}
$\psi$-{\em Appell} or generalized {\em Appell polynomials}
$\left\{A_{n}(x)\right\}_{n \geq 0}$ are equivalently
characterized via their $\psi$- exponential generating function
\begin{equation}\label{psi-A-3}
\sum_{n \geq
0}z^{n}\frac{A_{n}(x)}{n_{\psi}!}=A(z)\exp_{\psi}\left\{xz\right\},
\end{equation}
where $A(z)$ is a formal series with constant term different from
zero - here normalized to one.

The $\psi$- exponential function of $\psi$-Appell-Ward numbers
$A_{n}=A_{n}(0)$ is
\begin{equation} \label{psi-A-4}
\sum_{n \geq 0}z^{n}\frac{A_{n}}{n_{\psi}!}=A(z).
\end{equation}
Naturally $\psi$- {\em Appell} $\left\{A_{n}(x)\right\}_{n \geq
0}$ satisfy the $\psi$- {\em difference} equation
\begin{equation}\label{psi-A-5}
QA_{n}(x)=n_{\psi}x^{n-1};\;\;\;'n\geq 0,
\end{equation}
because
$QA_{n}(x)=QS^{-1}x^{n}=Q(\partial_{\psi}/Q)x^{n}=\partial_{\psi}x^{n}=n_{\psi}x^{n-1}\;
;n \geq 0$. Therefore they play the same role in
$Q(\partial_{\psi})$- {\em difference} calculus as Bernoulli
polynomials do in standard difference calculus or $\psi$-{\em
Bernoulli-Ward polynomials} (see Theorem 16.1 in \cite{1} and
consult also \cite{8}) in $\psi$-{\em difference} calculus due to
the following: The central problem of the $Q(\partial_{\psi})$ -
{\em difference calculus} is:
$$ Q(\partial_{\psi})f=\varphi\;\;\;\;\;\;\;\;\;\;\varphi=?,$$
where $f, \varphi$ - are for example formal series or polynomials.

The idea of finding solutions is the $\psi$-{\em Finite Operator
Calculus} \cite{2,3,4,5} standard. As we know (Proposition 2.1,
see \cite{2,3}) any $\psi$- delta operator $Q$ is of the form
$Q(\partial_{\psi})=\partial_{\psi}S$ where$S\in \Sigma_{\psi}$.
Let $Q(\partial_{\psi})=\sum_{k \geq
1}\frac{q_{k}}{k_{\psi}!}\partial_{\psi}^{k},\;\;q_{1}\neq 0$.
Consider then $Q(\partial_{\psi}=\partial_{\psi}S)$ with
$S=\sum_{k \geq
0}\frac{q_{k+1}}{(k+1)_{\psi}!}\partial_{\psi}^{k}\equiv \sum_{k
\geq
0}\frac{s_{k}}{k_{\psi}!}\partial_{\psi}^{k};\;\;\;s_{0}=q_{1}\neq
0$. We have for $S^{-1}\equiv\hat{A}$ - call it: $\psi$- {\em
Appell operator} - the obvious expression
$$\hat{A}\equiv S^{-1}=\frac{\partial_{\psi}}{Q_{\psi}}=\sum_{n\geq 0}\frac{A_{n}}{n_{\psi}!}\partial_{\psi}^{n}.$$
Now multiply the equation $Q(\partial_{\psi}f=\varphi)$ by
$\hat{A}\equiv \sum_{n\geq
0}\frac{A_{n}}{n_{\psi}!}\partial_{\psi}^{n}$ thus getting
\begin{equation}\label{psi_A_6}
\partial_{\psi}f=\sum_{n\geq
0}\frac{A_{n}}{n_{\psi}!}\varphi^{(n)},\;\;\;\varphi^{(n)}=\partial_{\psi}\varphi^{(n-1)}.
\end{equation}
The solution then reads:
\begin{equation}\label{psi-A-7}
f(x)=\sum_{n \geq
1}\frac{A_{n}}{n_{\psi}!}\varphi^{(n-1)}(x)+\int_{\psi}\varphi(x)+p(x),
\end{equation}
where $p$ is "$Q(\partial_{\psi})$- periodic" i.e.
$Q(\partial_{\psi})p=0$. Compare with  \cite{8} for "$+_{\psi}1$-
periodic" i.e. $p(x+_{\psi}1)=p(x)$ i.e. $\Delta_{\psi}p=0$. Here
 the relevant  $\psi$ - integration $\int_{\psi}\varphi(x)$ is defined as in
\cite{2}. We recall it in brief. Let us introduce the following
representation for $\partial_{\psi}$ "difference-ization"
$$ \partial_{\psi}=\hat{n}_{\psi}\partial_{0}\; ;\;\;\;\hat{n}_{\psi}x^{n-1}=n_{\psi}x^{n-1};\;\;n\geq 1,$$
where $\partial_{0}x^{n}=x^{n-1}$ i.e. $\partial_{0}$ is the $q=0$
Jackson derivative. $\partial_{0}$ is identical with divided
difference operator. Then we define the linear mapping
$\int_{\psi}$ accordingly:
$$ \int_{\psi}x^{n}=\left( \hat{x}\frac{1}{\hat{n}_{\psi}}\right)x^{n}=\frac{1}{(n+1)_{\psi}}x^{n+1};\;\;\;n\geq 0$$
where of course  $\partial_{\psi}\circ \int_{\psi}=id$.
\section{Examples}
\renewcommand{\labelenumi}{(\alph{enumi})}
\begin{enumerate}
\item The case of $\psi$- {\em Bernoulli-Ward} polynomials and
$\Delta_{\psi}$- {\em difference calculus} was considered in
detail in \cite{8} following \cite{1}. \item Specification of (a)
to the Gauss and Heine originating $q$-umbral calculus case
\cite{1,2,3,4,5} was already presented in \cite{8}. \item
Specification of (a) to the Lucas originating FFOC - case was also
presented in \cite{8} (here: FFOC={\bf F}ibonomial {\bf F}inite
{\bf O}perator {\bf C}alculus), see example 2.1 in \cite{5}).
Recall: the {\em Fibonomial coefficients} (known to Lucas)
($F_{n}$- {\em Fibonacci numbers}) are defined as
$$\binom{n}{k}_{F}=\frac{F_{n}!}{F_{k}!F_{n-k}!}=\binom{n}{n-k}_{F},$$
where in up-side down notation: $n_{F}\equiv F_{n}\neq
0$,\\$n_{F}!=n_{F}(n-1)_{F}(n-2)_{F}(n-3)_{F}\ldots
2_{F}1_{F};\;\;0_{F}!=1$;\\$n_{F}^{\underline{k}}=n_{F}(n-1)_{F}\ldots
(n-k+1)_{F};\;\;\;\binom{n}{k}_{F}\equiv
\frac{n^{\underline{k}}_{F}}{k_{F}!}$. We shall call the
corresponding linear difference operator
$\partial_{F};\;\;\partial_{F}x^{n}=n_{F}x^{n-1};\;\;n\geq 0$ the
$F$-derivative. Then in conformity with \cite{1} and with notation
as in \cite{2}-\cite{6} one has:
$$E^{a}(\partial_{F})=\sum_{n\geq 0}\frac{a^{n}}{n_{F}!}\partial_{F}^{n}$$
for the corresponding generalized translation operator
$E^{a}(\partial_{F})$. The $\psi$- integration becomes now still
not explored $F$- integration and we arrive at the$F$- Bernoulli
polynomials unknown till now.\\ {\bf Note:} recently a
combinatorial interpretation of Fibonomial coefficient has been
found \cite{6,7} by the present author.
 \item The other examples of $Q(\partial_{\psi})$- {\em difference calculus} -
expected naturally to be of primary importance in applications
(for inspiration see [1] and [9]- functional formulation) are
provided by the possible use of such $\psi$- Appell polynomials
as:

\begin{itemize}

\item

$\psi$-Hermite polynomials $\left\{ H_{n,\psi}\right\}_{n\geq 0}$:
$$ H_{n,\psi}(x)=\left[\sum_{k\geq 0}\left(-\frac{1}{2}\right)^{k}\frac{\partial_{\psi}^{2k}}{k_{\psi}!}\right]x^{n}\;\;\;n\geq 0;$$

\item

$\psi$ - Laguerre polynomials $\left\{ L_{n,\psi}\right\}_{n\geq
0}$ \cite{3}:

\begin{multline*}
L_{n,\psi}(x) = \frac{{n_{q}} }{{n}}\hat {x}_{\psi}
\left[{\frac{{1}}{{\partial _{\psi} - 1}}}\right]^{-n}x^{n-1} =
\frac{{n_{\psi}} }{{n}}\hat {x}_{\psi} \left( {\partial _{\psi} -
1} \right)^{n}x^{n-1} = \\ = \frac{{n_{\psi}} }{{n}}\hat
{x}_{\psi} \sum\limits_{k = 1}^{n} \left( { - 1} \right)^{k}\left(
{{\begin{array}{*{20}c}
 {n} \hfill \\
 {k} \hfill \\
\end{array}} } \right) \partial _{\psi} ^{n - k}x^{n - 1} =\\
=\frac{{n_{\psi}} }{{n}}\sum\limits_{k = 1}^{n} {} \left( { - 1}
\right)^{k}\left( {{\begin{array}{*{20}c}
 {n} \hfill \\
 {k} \hfill \\
\end{array}} } \right) \left( {n - 1} \right)_{\psi} ^{\underline {n - k}
}\frac{{k}}{{k_{\psi}} }x^{k} .
\end{multline*}

For $q=1$ in $q$-extended case [9] one recovers the known formula
:

\begin{center}
$L_{n,q=1}(x) = \sum\limits_{k = 1}^{n} {} \left( { - 1}
\right)^{k}\frac{{n_{q} !}}{{k_{q} !}}\left(
{{\begin{array}{*{20}c}
 {n - 1} \hfill \\
 {k - 1} \hfill \\
\end{array}} } \right)x^{k}$.
\end{center}

\end{itemize}
\end{enumerate}


\begin{thebibliography}{10}

\bibitem{1}
M. Ward: {\em A calculus of sequences}, Amer.J.Math. Vol.58,
1936,pp.255-266

\bibitem{2}
A. K. Kwa\'sniewski: {\em Main Theorems of Extended Finite
Operator Calculus}, Integral Transforms and Special Functions Vol
14, No 6, (2003), pp.499-516

\bibitem{3}
A. K. Kwa\'sniewski: {\em Towards $\psi$-extension of Finite
Operator Calculus of Rota}, Rep.Math.Phys. 48 No3 (2001)
pp.305-342 ArXiv:math.CO/0402078 2004

\bibitem{4}
A. K. Kwa\'sniewski: {\em On Extended Finite Operator Calculus of
Rota and Quantum Groups}, Integral Transforms and Special
Functions Vol 2, No 4, (2001), pp.333-340

\bibitem{5}
A. K. Kwa\'sniewski: {\em On Characterizations of Sheffer
$\psi$-polynomials and Related Propositions of the Calculus of
Sequences}, Bulletin de la Soc. des Sciences et des Lettres de
Lodz, 52, Ser.Rech.Deform.36(2002)pp.45-65,  ArXiv:math.CO/0312397

\bibitem{6}
A. K. Kwa\'sniewski: {\em Information on Combinatorial
Interpretation of Fibonomial coefficients}, Bulletin de la Soc.
des Sciences et des Lettres de Lodz, 53,
Ser.Rech.Deform.42(2003)pp.39-41

\bibitem{7}
A.K.Kwa\'sniewski: {\em  More on combinatorial interpretation of
the fibonomial coefficients}, ArXiv: math.CO/0402344 v1  2004

\bibitem{8}
A. K. Kwa\'sniewski: {\em A note on $\psi$-basic Bernoulli-Ward
Polynomials and Their Specifications}, ArXiv: math.CO/0405577 30
May 2004

\bibitem{9}
S. M. Roman: {\em The umbral calculus\/}, Academic Press, New York
1984.

\end{thebibliography}
 \end{document}